\documentclass[11pt]{cmslatex}

\usepackage[margin=1in, top=1.5in]{geometry}
\usepackage{amssymb}
\usepackage{amsmath}
\usepackage{bbm}
\usepackage{graphicx}

\renewcommand{\theequation}{\arabic{section}.\arabic{equation}}

            \newtheorem{thm}{Theorem}[section]
          
          \newtheorem{lem}[thm]{Lemma}

          \newtheorem{defn}[thm]{Definition}

\newcommand{\id}{\mathbbm{1}}
\renewcommand{\grad}{\nabla{}}
\renewcommand{\div}{\mbox{\rm div}}

\newcommand{\e}{\varepsilon}

\title{Approximate solutions of scalar conservation laws}
 \author{Misha Perepelitsa \thanks{ misha@math.uh.edu,
University of Houston,
PGH 631,
4800 Calhoun Rd.,
Houston, TX.
USA}
        }

\begin{document}

\pagestyle{myheadings} \markboth{Approximate solutions of conservation laws}{Misha Perepelitsa}

\maketitle
\begin{abstract}
We study compactness properties of time-discrete and continuous time BGK-type schemes for scalar conservation laws, in which microscopic interactions occur only when the state of
a system deviates significantly from an equilibrium distribution. The threshold deviation, $\e,$ is a parameter of the problem. In the vanishing relaxation time limit we obtain solutions of a conservation law in which flux is pointwisely close (of order $\e$) to the flux of the original equation and derive several other properties of such solutions, including an example of approximate solution to a shock for Burger's equation.
\end{abstract}

\begin{keywords}
Scalar conservation laws, kinetic models, shock waves.
\end{keywords}

 \begin{AMS}
35L02, 35L65,  58J45.
\end{AMS}

\begin{section}{Introduction}
Several important classes of PDEs can be formally represented by using kinetic formalism by extending the space of independent variables. In fact, many PDEs, such as the Euler equations of gas dynamics can be derived from the kinetic equations of microscopic description of motion. A kinetic representation of equations of the type
\begin{equation}
\label{eq:system}
\partial_t U{}+{}\div_x F(U){}={}0,
\end{equation}
where $(x,t)\in\mathbb{R}^{d+1},$ $U\;:\;\mathbb{R}^{d+1}\to \mathbb{R}^d,$ and $F(U)\;:\;\mathbb{R}^d\to \mathbb{R}^{d\times d},$ uses a microscopic variable $v$ and a set of equilibria $\mathcal{M}_0$ that consist of functions $f(U,v)$ such that for some functions $\kappa(v)\in\mathbb{R}^d$ and $a(v)\in\mathbb{R}^{d},$ such that for any $U\in\mathbb{R}^d,$
\[
U{}={}\int \kappa(v)f(U,v)\,dv,\quad F(U){}={}\int \kappa(v)\otimes a(v) f(U,v)\,dv.
\]
One attempts to solve the kinetic equation for unknowns $(f,Q):$
\begin{equation}
\label{eq:system_kinetic}
\partial_t f{}+a(v)\cdot\grad_x f{}={}Q,
\end{equation}
such that for all $(x,t)$
\[
f(x,t,\cdot)\in \mathcal{M}_0,\quad \int \kappa(v) Q(x,t,v)\,dv{}={}0.
\]
Often more information is available: the set of equilibria $\mathcal{M}_0$ can be characterized as a set of constrained  minimizers of an entropy (convex function) $S(f):$
\[
\min \{ S(f)\;:\; \int \kappa(v)f(v)\,dv{}={}U=const.\}.
\]
For such $f, $ $\int S(f(U,v))\,dv{}={}s(U)$ is an entropy function for \eqref{eq:system}:
\[
\partial_t s{}+{}\div_x q\leq0,
\]
where $q=q(U)$ is the entropy flux. We will use notation $\Pi_f\in \mathcal{M}_0$ for the above minimizer with $\kappa$--moments: $\int\kappa(v)\Pi_f(v)\,dv{}={}\int \kappa(v)f(v)\,dv.$

The results \cite{LPT, LPT2} show that admissible weak solutions of scalar conservation laws in multi dimensions, and equations of isentropic gas dynamics is dimension one,
can be represented by kinetic densities that verify \eqref{eq:system_kinetic} with $Q$ in the form $-\partial_v m,$ for scalar conservation laws, or $\partial^2_v m,$ for isentropic gas dynamics, where $m$ is non-negative Radon's measure on space $(x,t,v).$

One approach to  solving  problem \eqref{eq:system_kinetic} consists in relaxing the requirement that $f(x,t,\cdot)\in \mathcal{M}_0$ and selecting a particular ``interaction potential" $Q.$ For example, in a BGK-type model 
\begin{equation}
\label{eq:systems_BGK}
\partial_t f{}+{}a(v)\cdot \grad_x f{}={}Q_h(f){}={}\frac{\Pi_f-f}{h},
\end{equation}
where $h>0$ is ``relaxation" time scale. When this equation is paired with  $\partial_fS(f),$ one obtains an entropy inequality
\begin{equation}
\label{eq:system_entropy}
\int_0^T\int_{\mathbb{R}^d_x}\langle \partial^2_f S(f)(\Pi_f-f),(\Pi_f-f)\rangle\,dxdt\leq h\left(\int_{\mathbb{R}^d_x} \int S(f(x,0,v))-S(f(x,T,v))\,dvdx\right).  
\end{equation}
It shows that $\int_{\mathbb{R}^d_x} \int S(f(x,t,v))\,dvdx$ doesn't grow in time, and on average (in $(x,t)$), $f(x,t,\cdot)$ is close to the equilibrium density $\Pi_f(\cdot).$

For many equations of interest, such as scalar conservation laws, isentropic gas dynamics, the Euler equation of gas dynamics the theory of BGK-type models is well--developed with results establishing existence/uniqueness and stability of solutions in appropriate spaces, see for example the monograph \cite{Perthame} and references therein.

One naturally interested in trying to recover solutions of \eqref{eq:system} by passing to the limit of zero relaxation time $h$ in a family of solutions $f_h$ of \eqref{eq:systems_BGK}, with some fixed initial data $f(x,0,v).$ Save for scalar conservation laws and one-dimensional isentropic gas dynamics, \cite{Brenier, GM, BB}, current techniques based on averaging lemmas do not allow to conclude that $f_h$ accumulates on  a solution of \eqref{eq:system}, for the reason the entropy estimate \eqref{eq:system_entropy} is not sufficient ``to control"  the right-hand side $Q_h$ in \eqref{eq:systems_BGK}.

In \cite{Perepelitsa_isentropic} the model \eqref{eq:systems_BGK} for multidimensional, isentropic gas dynamics  was further relaxed in an attempt to obtain solutions $f$ that are pointwisely close to the equilibrium density $\Pi_f.$ The closeness is measured by a relative change of the entropy:
\[
D(f(x,t,\cdot)){}={}\int S(f(x,t,\cdot))- S(\Pi_{f(x,t,\cdot)}(v))\,dv/ \int S(\Pi_{f(x,t,\cdot)}(v))\,dv{}\leq{}\e,
\]
where $\e$ is a non-dimensional parameter.

The following  BGK-type model was considered in \cite{Perepelitsa_isentropic}
\begin{equation}
\label{eq:systems_BGK_new}
\partial_t f{}+{}v\cdot \grad_x f{}={}Q_h(f){}={}\frac{\Pi_f-f}{h}\id_{\{D(f(x,t,\cdot))>\e\}}(x,t),
\end{equation}
and was shown that a family $\{f_h\}$ of its solutions accumulates on a kinetic function $f$ that verifies equation
\[
\partial_t f{}+{}v\cdot\grad_x f{}={}Q,
\]
with $Q$ -- signed Radon's measure. Such $f$ is close to the equilibrium density: for a.e. $(x,t),$ $D(f(x,t,\cdot))\leq \e.$ Moreover, the moments $(\rho,\rho u){}={}\int (1,v)f(x,t,v)\,dv$ solve the system of approximately isentropic Euler equations with the pressure function $p(x,t)$ that is comparable to the isentropic pressure $\kappa\rho^\gamma$ as:
\[
\kappa \rho^\gamma(x,t)\leq p(x,t)\leq (1+O(\e))\kappa\rho^\gamma(x,t).
\]

Since the concept of an approximate solutions of this type is rather general, in the present work, we apply it to scalar conservation laws and consider solutions in some details.

\subsection{Scalar conservation laws}
In the half-plane $\mathbb{R}^+_{x,t}{}={}(\mathbb{R}^d\times \mathbb{R}^+),$ $d\geq1,$ we consider Cauchy problem for a scalar conservation law
\begin{equation}
\label{eq:SCL} 
\partial_t u{}+{}\div_x A(u){}={}0,\quad (x,t)\in\mathbb{R}^+_{x,t},
\end{equation}
 with the initial conditions
\[
u(x,0){}={}u_0(x),\quad x\in\mathbb{R}^d.
\]
In the above equation $A(u)\;:\;\mathbb{R}\to\mathbb{R}^d{}={}(A^1(u),...,A^d(u))$ is a flux function, which we always assume to be at least locally Lipschitz continuous. 

Let $\eta(v)$ be a convex (entropy) function on $\mathbb{R}_v$ with its flux defined as
\[
q(v){}={}\int_0^v \eta'(v) A'(v)\,dv.
\]

\begin{defn} A bounded measurable function $u\in L^\infty(\mathbb{R}^+_{x,t})$ is an admissible solution of the Cauchy problem for \eqref{eq:SCL} with $u_0\in L^\infty(\mathbb{R}^d),$
if for any $(\eta,q)$ -- entropy/entropy flux pairs and all non-negative Lipschitz continuous test function $\psi$ on $\mathbb{R}^d_x\times[0,+\infty),$ with compact support it holds:
\[
\iint_{\mathbb{R}^+_{x,t}} (\partial_t\psi \eta(u)+ \grad_x\psi\cdot q(u))\,dxdt{}+{}\int_{\mathbb{R}^d_x}\psi u_0\,dx \geq0.
\]
\end{defn}

The theory of admissible solutions is well-developed and is briefly summarized in the next theorem. 
\begin{thm} For any $u_0\in L^\infty(\mathbb{R}^d),$ there exists a unique admissible weak solution of \eqref{eq:SCL} and \[u\in C([0,\infty);L^1_{loc}(\mathbb{R}^d))\cap L^\infty(\mathbb{R}^+_{x,t}).\]
\end{thm}
The exposition of the theorem can be found, for example, in the monograph \cite{Constantine}, Ch. IV, which also  references the important contributions to this theory.
 The result of Lions-Perthame-Tadmor\cite{LPT} establishes an equivalence between  the class of  entropy solutions $u\in L^\infty(\mathbb{R}^+_{x,t})$ of \eqref{eq:SCL}, and the class of solutions in the form
\begin{equation}
\label{def:projection}
f(x,t,v){}={}
\left\{
\begin{array}{rl}
\id_{[0,u(x,t)]}(v), & u(x,t)\geq 0,\\
-\id_{[u(x,t),0]}(v), & u(x,t)<0.
\end{array}
\right.
\end{equation}
of the equation
\begin{equation}
\label{eq:kineticLPT}
\partial_t f{}+{}A'(v)\cdot\grad_x f{}={}-\partial_v m,
\end{equation}
with $m$ -- a non-negative Radon measure on $\mathbb{R}^+_{x,t,v}$  The equivalence is defined through \eqref{def:projection}. The set of equilibrium densities $\mathcal{M}_0,$ mentioned above, consists of functions of this type.

One way to construct admissible solutions is through  the kinetic time-discrete scheme of Brenier\cite{Brenier} and Giga-Miyakawa\cite{GM}. The scheme is  defined in the following way. Given $u_0(x),$ set the kinetic density
\begin{equation}
\label{projection}
f(x,0,v){}={}
\left\{
\begin{array}{rl}
\id_{[0,u_0(x)]}(v), & u_0(x)\geq 0,\\
-\id_{[u_0(x),0]}(v), & u_0(x)<0.
\end{array}
\right.,\quad \int f(x,0,v)\,dv{}={}u_0(x).
\end{equation}
Given time step $h>0$ define function $f_h(x,t,v)$ recursively on intervals $(nh, (n+1)h)$ by first evolving it by transport:
\begin{equation}
\label{step:transportB}
f_h(x,t,v){}={}f_h(x-(t-nh)A'(v),nh,v),\quad t\in(nh,(n+1)h).
\end{equation}
Then projecting it to an equilibrium kinetic density, preserving the zero moment in $v:$
\begin{equation}
\label{step:collapseB}
f_h(x,(n+1)h,v){}={}
\left\{\begin{array}{rl}
\id_{[0,\tilde{u}(x)]}(v), & u_0(x)\geq 0,\\
-\id_{[\tilde{u}(x),0]}(v), & u_0(x)<0,
\end{array}
\right.
\end{equation}
where
\[
\tilde{u}(x){}={}\int f_h(x-hA'(v),nh,v)\,dv.
\]

It was shown in \cite{Brenier} that if $A(u)$ is Lipschitz continuous, for any $u_0\in L^1(\mathbb{R}^d),$
the sequence of functions
\[
u_h(x,t){}={}\int f_h(x,t,v)\,dv
\]
converges in topology of $C([0,T];L^1(\mathbb{R}^d))$ for any $T>0,$ to a function $u(x,t),$
which is a unique entropy solution of \eqref{eq:SCL}, and there is a non-expansive semi-group on $S_t$ on $L^1(\mathbb{R}^d)$ such that
\[
u(x,t){}={}S_t[u_0](x).
\]
Moreover, the semi-group $S_t $ is order preserving: if $u_0(x)\leq u_1(x)$ $a.e.\,x,$ then
for any $t>0,$
\[
S_t[u_0](x)\leq S_t[u_1](x),\quad a.e.\,x.
\]

The change of an entropy in this scheme is controlled in the transport step by equation
\[
\partial_t \left(\int \eta'(v)f_h\,dv\right){}+{}\div_x \left( \int \eta'(v)A'(v) f_h\,dv\right){}={}0,\quad \mbox{in } \mathcal{D}'(\mathbb{R}^d_x\times  (nh,(n+1)h)),\forall n,
\]
and in the collapse step by an inequality 
\begin{equation}
\label{0:1}
\int \eta'(v)(f_h(x,(n+1)h,v)-f_h(x-hA'(v), nh,v))\,dv\geq0,\quad \forall x,
\end{equation}
where $\eta(v)$ is a convex function. Combined,  
\[
\partial_t \left(\int \eta'(v)f_h\,dv\right){}+{}\div_x \left( \int \eta'(v)A'(v) f_h\,dv\right){}\leq{}0,\quad \mbox{in }\mathcal{D}'(\mathbb{R}^+_{x,t}).
\]

 The inequality \eqref{0:1} follows from a convexity lemma, \cite{Brenier}, p. 1017.  For a reference we recall the lemma in the Appendix, lemma \ref{lemma:Brenier}.

\subsection{Main results} In this paper we consider a modification of the above scheme in which the projection step is
performed only for $x$ where the decrease of an entropy \eqref{0:1} is large, as described below. We select quadratic function $\eta(v){}={}v^2/2$ as a designated entropy and measure its decrease in the ``projection step''  \eqref{step:collapseB} by the ratio
\[
D(f(\cdot)){}={} \int v(f(v)-\Pi_{f(\cdot)}(v))\,dv/ \int v \Pi_{f(\cdot)}(v))\,dv\geq0.
\]
Function $D(f)$ is defined by zero if $\int v \Pi_{f(\cdot)}(v))\,dv{}={}0.$

Fix $\e>0.$  Given time step $h>0$ define function $f_h(x,t,v)$ recursively on intervals $(nh, (n+1)h)$ by first evolving it by transport:
\begin{equation}
\label{step:transport}
f_h(x,t,v){}={}f_h(x-(t-nh)A'(v),nh,v),\quad t\in(nh,(n+1)h).
\end{equation}
Then, for points $x$ such that
\begin{equation}
\label{step:small}
D(f_h(x-hA'(v),nh,\cdot))>\e,
\end{equation}
we  project $f_h(x-hA'(v),nh,\cdot)$ to an equilibrium kinetic density 
\begin{equation}
\label{step:collapse}
f_h(x,(n+1)h,v){}={}\Pi_{f_h(x-hA'(\cdot),nh,\cdot)}(v)
\end{equation}

In the following results we always assume that $A\in C^2(\mathbb{R}_v)^d$ and verifies a non-degeneracy condition:
\begin{equation}
\label{main:nondegeneracy}
\forall\, \sigma \in \mathbb{S}^{d-1}\quad \forall \xi\in\mathbb{R}\quad |\{ v\in[-\|u_0\|_{L^\infty},\|u_0\|_{L^\infty}]\;:\; A'(v)\cdot\sigma{}={}\xi\}|{}={}0,
\end{equation}
where $\mathbb{S}^{d-1}$ is the unit sphere in $\mathbb{R}^d.$
\begin{thm}
\label{th:1}
Let $u_0\in L^\infty(\mathbb{R}^d_x).$ For the sequence $\{f_h\}$ defined above there is a subsequence (still labeled by $h$) such that as $h\to0,$
\begin{enumerate}
\item $f_h\to f$ *--weakly in $L^\infty_{loc}(\mathbb{R}^+_{x,t,v})$ to some function $f$ with values in $[-1,1];$
\item for all test functions $\phi\in C(\mathbb{R}_v),$
\[
\left\{ \int \phi f_h\,dv \right\}\quad  \mbox{compact in $L^p_{loc}(\mathbb{R}^+_{x,t})$}
\]
for all $1\leq p<+\infty;$
\item for every entropy $\eta,$
\[
\partial_t \int \eta' f\,dv{}+{}\div_x \int \eta' A'f\,dv\leq 0,\quad \mathcal{D}'(\mathbb{R}^+_{x,t});
\]
\item there is a non-negative Radon measure on $\mathbb{R}^+_{x,t,v}$ such that
\[
\partial_t f{}+{}A'\cdot\grad_x f{}={}-\partial_v m.
\]
Moreover, for any compact $K\subset \mathbb{R}^+_{x,t,v},$
$m[K]$ is bounded independently of $\e;$
\item if a.e. $x,$ 
\[
u_0(x)\geq \bar{u}>0\,(\mbox{or } u_0(x)\leq -\bar{u}<0),
\] then the right-hand side of the above kinetic equation, $\partial_v m,$ is itself a signed Radon's measure (signed Borel's measure with locally finite mass).
\end{enumerate}
\end{thm}
Note that there is gain of a derivative in the right-hand side of the kinetic equation, as compared with the kinetic equation of corresponding to an admissible solution of \eqref{eq:SCL}.

The proof uses entropy inequalities, optimal transport maps between to kinetic densities and a compactness theorem of G\'{e}rard\cite{Gerard}. The proof of this the subsequent theorems are given the next section.

In the next theorem we show that $u(x,t){}={}\int f(x,t,v)\,dv$ from theorem \ref{th:1} solves an equation which is approximately close to \eqref{eq:SCL}, and can be represented by a measure-valued solution of \eqref{eq:SCL}, 
which is close to a Dirac's mass concentrated on $u(x,t):$ $\delta(v-u(x,t)).$

\begin{thm}
\label{th:2}
Let $f(x,t,v)$ be a function verifying parts (3), (4) of theorem \ref{th:1} and  such that a.e. $(x,t),$
\begin{equation}
\label{cond:small}
D(f(x,t,\cdot))\leq \e.
\end{equation}
Let $u(x,t){}={}\int f(x,t,v)\,dv.$
Then,
\\
\noindent (i) $u$ is an approximate solution of \eqref{eq:SCL} in a sense that for any strictly convex entropy $\eta,$ (normalized by $\eta(0){}={}0$) and its flux $q,$
\[
\partial_t \left(\eta(u)(1+O_\eta(\e))\right){}+{}\div_x\left( \eta(u)(\frac{q(u)}{\eta(u)}+O_q(\e))\right)\leq 0,\quad
\mbox{in } \mathcal{D}'(\mathbb{R}^+_{x,t}),
\]
where $O_\eta(\e), O_q(\e)$ stand for    functions of $(x,t)$ such that
\[
\|O_\eta(\e)\|_{L^\infty(\mathbb{R}^+_{x,t})},\,\|O_q(\e)\|_{L^\infty(\mathbb{R}^+_{x,t})} {}\leq{}C\e.
\]
In particular,
\[
\partial_t u{}+{}\div_x\left( u(A(u)/u+O(\e))\right){}={}0,\quad
\mbox{in } \mathcal{D}'(\mathbb{R}^+_{x,t});
\]
\noindent (ii) there is a parametrized signed, unit mass measure $\mu_{x,t}{}={}\delta(v-u(x,t)){}+{}m_{x,t}$ such that 
\[
\mbox{\rm mass}(m_{x,t}){}={}O(\e),\quad \mbox{a.e. }(x,t),
\]
and $\mu_{x,t}$ is a measure-valued solution of \eqref{eq:SCL}:
\[
\partial_t \int v\,d\mu_{x,t}{}+{}\div_x\int A(v)\,d\mu_{x,t}{}={}0,\quad \mbox{in }\mathcal{D}'(\mathbb{R}^+_{x,t}).
\]

\end{thm}

It is not clear if condition \eqref{cond:small} holds in the limit of the above discrete time scheme, but it  can be obtained in the zero relaxation limit of a BGK-type scheme that we describe next.

Consider the problem
\begin{equation}
\label{eq:2.1}
\partial_t f{}+{}A'(v)\cdot\grad_x f{}={}\frac{(\Pi_{f(x,t,\cdot)}(v)-f(x,t,v))}{h}\id_{\{D(f(x,t,\cdot))>\e\}}(x,t),\quad \mathcal{D}'(\mathbb{R}^+_{x,t,v}),
\end{equation}
\begin{equation}
\label{eq:2.2}
f(x,0,v){}={}f_0(x,v).
\end{equation}
For every $h>0$ and $f_0$ from an appropriate space, the problem has a weak solutions, see theorem \ref{th:last} in the Appendix. For a family $\{f_h\}$ of solutions of 
\eqref{eq:2.1}, \eqref{eq:2.2} we prove the following
\begin{thm} 
\label{th:4}
Let $u_0\in L^\infty(\mathbb{R}^d_x)\cap L^1(\mathbb{R}^d_x),$ and $f_0$ be the corresponding kinetic density. For the family $\{f_h\}$ of weak solutions of \eqref{eq:2.1}, \eqref{eq:2.2}, the following holds.
\begin{enumerate}
\item $f_h\to f$ *--weakly in $L^\infty(\mathbb{R}^+_{x,t,v})$ to some function $f$ with values in $[-1,1],$ and 
\[
f\in L^\infty(\mathbb{R}^+_{x,t,v})\cap L^\infty((0,+\infty); L^1(\mathbb{R}^{d+1}_{x,v}));
\]
\item for all test functions $\phi\in C(\mathbb{R}_v),$
\[
\left\{ \int \phi f_h\,dv \right\}\quad  \mbox{compact in $L^p_{loc}(\mathbb{R}^+_{x,t})$}
\]
for all $1\leq p<+\infty;$
\item for every entropy $\eta,$
\[
\partial_t \int \eta' f\,dv{}+{}\div_x \int \eta' A'f\,dv\leq 0,\quad \mathcal{D}'(\mathbb{R}^+_{x,t});
\]
\item there is a non-negative Radon measure on $\mathbb{R}^+_{x,t,v}$ such that
\[
\partial_t f{}+{}A'\cdot\grad_x f{}={}-\partial_v m.
\]
Moreover, for any compact $K\subset \mathbb{R}^+_{x,t,v},$
$m[K]$ is bounded independently of $\e;$
\item for  a.e. $(x,t),$
\[
D(f(x,t,\cdot))\leq \e;
\]
\item if a.e. $x,$ 
\[
u_0(x)\geq \bar{u}>0\,(\mbox{or } u_0(x)\leq -\bar{u}<0),
\] then the right-hand side of the above kinetic equation, $\partial_v m,$ is itself a signed Radon's measure (signed Borel's measure with locally finite mass).
\end{enumerate}
\end{thm}

Lastly, we give an example of a self-similar approximate solution of  Burger's equation, corresponding to a stationary shock.

\end{section}

\begin{section}{Proofs}
\proof[Proof of theorem \ref{th:1}] Without loss of generality we assume that $u_0(x)$ is non-negative for a.e. $x$ and denote $M{}={}\|u_0\|_{L^\infty}.$ The family of kinetic densities $\{f_h\}$ described in the introduction is well-defined on $\mathbb{R}^+_{x,t,v}$ and has the following properties that are easily  verified.
\begin{enumerate}
\item For all $(x,t,v),$ $f_h(x,t,v)\in[0,1];$
\item for all $(x,t),$  $ \supp f_h(x,t,\cdot)\subset[0,M];$
\item $f_h$ is a solution in a sense of distributions of equation
\begin{equation}
\label{eq:app_kinetic}
\partial_t f_h{}+{}A'\cdot\grad_x f_h{}={}\sum_{n=1}^{\infty }(f_h(x,t,v)-f_h,x,t-,v)   \id_{\{ D(f_h(x,t-,\cdot)>\e\} }(x,t)\delta(t-nh).
\end{equation}
\end{enumerate}
In the next lemma we list the energy (entropy) estimates. Recall that $f_h(x,nh,\cdot)$ is an equilibrium density obtained from $f_h(x,nh-,\cdot).$  Using lemma \ref{lemma:ap1} from Appendix we denote by $\tau_h(x,nh,v)$ the transport map that moves $ f_h(x,nh,v)\,dv$ to  $f_h(x,nh-,v)\,dv.$

\begin{lem} Let $s(x,t)$ be a smooth cut-off function for set $B_r(x)\times(r^{-1},r)\subset \mathbb{R}^+_{x,t},$ i.e. $s\in C^\infty_c(\mathbb{R}^+_{x,t})$ and $s=1$ on $B_r(x)\times(r^{-1},r)\subset \mathbb{R}^+_{x,t}.$ There is $C_r>0,$ independent of $(\e,h)$ such that
\begin{equation}
\label{est:1}
\sum_{n=\lceil 1/(rh)\rceil}^{\lfloor r/h\rfloor } \int_{B_r}  \id_{\{D(f_h(x,nh-,\cdot)>\e\}}(x) \int |\tau_h(x,nh,v)-v| f_h(x,nh,v)\,dvdx{}\leq{}C_r.
\end{equation}
There is $C_{r,\e}>0$ independent of $h,$ such that
\begin{equation}
\label{est:2}
\sum_{n=\lceil 1/(rh)\rceil}^{\lfloor r/h\rfloor } \int_{B_r}  \id_{\{D(f_h(x,nh-,\cdot)>\e\}}(x)  \int v f_h(x,nh,v)\,dvdx{}\leq{}C_{r,\e}.
\end{equation}
If a.e. $x,$ $u_0(x)\geq \bar{u}>0,$ then $[0,\bar{u}]\subset \supp f_h(x,nh,\cdot),$ and it follows from \eqref{est:2} that there is $C_{r,\e,\bar{u}}>0$ such that 
\begin{equation}
\label{est:3}
\sum_{n=\lceil 1/(rh)\rceil}^{\lfloor r/h\rfloor }\left| B_r\cap \{D(f_h(x,nh-,\cdot)>\e\} \right|{}\leq{}C_{r,\e,\bar{u}}.
\end{equation}
\end{lem}
\proof The estimates are straightforward, obtained by multiplying equation \eqref{eq:app_kinetic}
by $vs(x,t)$ and integrating over $\mathbb{R}^+_{x,t,v}.$

From this lemma we get the $v$--moments estimates for the right-hand side of \eqref{eq:app_kinetic}.
\begin{lem}
Let $\psi\in C^\infty_c(\mathbb{R}_v).$ Then 
\[
\int \psi(\partial_t f_h + A'\cdot\grad_x f_h)\,dv
\]
is a signed Radon measure on $\mathbb{R}^+_{x,t}$ with local total mass bounded independently of $(\e,h).$ 

If for  a.e. $x,$ $u_0(x)\geq \bar{u}>0,$ then the right-hand side of \eqref{eq:app_kinetic} is a signed Radon measure with local total mass bounded independently of $h.$
\end{lem}
\proof The first statement follows by pairing equation \eqref{eq:app_kinetic} with $\psi(v)s(x,t),$ and using estimate \eqref{est:1}. The second statement follows from \eqref{est:3} and fact that $f_h$ is uniformly bounded in $(x,t,v),$ independently of $h.$

 Now, since 
\[
 \int \psi(\partial_t f_h + A'\cdot\grad_x f_h)\,dv{}={}
\partial_t  \int \psi f_h\,dv + div_x \int A' \psi f_h\,dv
\]
 is also bounded in $W^{-1,\infty}_{loc}(\mathbb{R}^+_{x,t}),$ the family  of moment 
$\int \psi(\partial_t f_h + A'\cdot\grad_x f_h)\,dv$ is compact in $H^{-1}_{loc}(\mathbb{R}^+_{x,t}).$ 

Theorem of G\'{e}rard, (theorem \ref{Gerard} from the Appendix), applies to assert that
\[
\int\phi f_h\,dv
\]
are compact in the strong topology of $L^2_{loc}(\mathbb{R}^+_{x,t}),$ and thus a.e., and 
in every $L^p_{loc}(\mathbb{R}^+_{x,t}),$ for $1\leq p<+\infty.$

\proof[Proof of theorem \ref{th:2}] Assume that $u_0(x)\geq0$ a.e. $x,$ and set $M{}={}\|u_0\|_{L^\infty}.$ Consider $f(x,t,\cdot)\in L^\infty((0,M)),$ that is defined for a.e. $(x,t)$ and set $u{}={}u(x,t){}={}\int f\,dv.$ For a strictly convex on $[0,M]$ entropy $\eta(v),$ and its flux $q(v),$ 
\[
\int \eta'(v) f\,dv{}={}\eta(u){}+{}\int \eta'(v)(f-\Pi_f)\,dv,\,\int \eta'(v)A'(v)f\,dv{}={}q(u){}+{}\int \eta'(v)A'(v)(f-\Pi_f)\,dv.
\]
By the smallness condition \eqref{main:nondegeneracy}, 
\[
\frac{\left|\int \eta'(v)(f-\Pi_f)\,dv\right|}{\eta(u)}{}\leq C\e,\,
\frac{\left|\int \eta'(v)A'(v)(f-\Pi_f)\,dv\right|}{\eta(u)}{}\leq C\e,
\]
with some $C>0$ independent of $(x,t),$ from which part (i) follows.

To prove part (ii), let let $A_0(v){}={}v$ and consider $d+1$ dimensional subspace of $L^2((0,M)),$
\[
H{}={}\left\{ \sum_{k=0}^d \alpha_k A'(v)\;:\; \alpha_k\in\mathbb{R},\,k=0..d\right\}.
\]
Consider the projection of $f-\Pi_f\in L^2((0,M))$ onto $H,$  and denote it by $f_{min}.$ It verifies  the equations:
\[
\int A'_i f_{min}\,dv{}={}\int A'_i(v)(f-\Pi_f)\,dv{}={}O(\e),\quad i=0,..,d.
\]
Since $f_{min}{}=\sum \alpha_k A_k'(v),$ for some $\alpha_k$'s. The coefficients are solutions of $d+1$ linear equations
\[
A\bar{\alpha}{}={}\bar{b},\quad A{}={}\{a_{ij}\},\, a_{ij}{}={}\int A'_iA'_j\,dv, \,\bar{\alpha}{}={}(\alpha_0,...,\alpha_d)^T,
\]
and $\bar{b}$ is the vector of right-hand side with each $b_i{}={}O(\e).$ 

The non-degenracy conditions \eqref{main:nondegeneracy} implies that the set $\{A'_0,A'_1,...,A'_d\}$ is linearly independent  on $[0,M],$ and so, matrix $A$ is non-singular, each $\alpha_i{}={}O(\e),$
and $\|f_{min}\|_{L^\infty},\,\|f_{min}'\|_{L^\infty}{}={}O(\e).$
A measure representing a solution \eqref{eq:SCL} is defined as
\[
\mu_{x,t}{}={}\delta(v-u){}-{}f_{min}'(v)\,dv+f_{min}(M)\delta(v-M).
\]

\proof [Proof of  theorem \ref{th:4}] Assume without the loss of generality that $u_0(x)\geq 0,$ a.e. $x.$ The proof of this theorem uses the same compactness argument as in theorem \ref{th:1}, based on the entropy estimates of the next lemma. As above, we denote by $\tau_h(x,t,v)$ the transport map from $\Pi_{f(x,t,\cdot)}(v)$ to $f(x,t,\cdot),$ from lemma \ref{lemma:ap1}.
\begin{lem} Let $B_r\subset{\mathbb{R}^d_x}$ be a ball of radius $r$ centered at the origin.
There is $C_r>0$ independent of $(\e,h),$ such that
\begin{equation}
\label{est:2.1}
\int_{1/r}^r \int_{B_r}\id_{\{D(f_h(x,t,\cdot)>\e\}}(x,t)\int |\tau_h-v|\Pi_{f_h(x,t,\cdot)}(v) \,dvdxdt {}\leq C_rh;
\end{equation}
there is $C_{r,\e}>0,$ independent of $h,$ such that
\begin{equation}
\label{est:2.2}
\int_{1/r}^r \int_{B_r}\id_{\{D(f_h(x,t,\cdot)>\e\}}(x,t)\int v\Pi_{f_h(x,t,\cdot)}(v) \,dvdxdt {}\leq C_{r,\e}h;
\end{equation}
If a.e. $x,$ $u_0(x)\geq \bar{u}>0,$ then $[0,\bar{u}]\subset \supp f_h(x,t,\cdot),$ a.e. $(x,t),$ and it follows from \eqref{est:2.2} that there is $C_{r,\e,\bar{u}}>0$ such that 
\begin{equation}
\label{est:2.3}
 |B_r\times(1/r,r)\cap \{D(f_h(x,t,\cdot))>\e\}|\,dt{}\leq{}C_{r,\e,\bar{u}}h.
\end{equation}
\end{lem}
The lemma is obtained by taking $v$--moment of the equation \eqref{eq:2.1}, and we omit it.
With this lemma all statements of the theorem are proved exactly as in the proof of theorem above, except part 5, which becomes evident once we observe \eqref{est:2.2}: for any $w(x,t)\in C^\infty_0(\mathbb{R}^+_{x,t}),$ 
\begin{multline*}
\int\int \left(\int v(f_h-\Pi_{f_h})\,dv\right)w(x,t)\id_{\{D(f_h(x,t,\cdot)>\e)\}}(x,t)\,dxdt{}={}\\
\int\int \left(\int v(f_h-\Pi_{f_h})\,dv\right)w(x,t)\id_{\{\int \Pi_{f_h}\,dv<h^{1/4}\}\cap\{D(f_h(x,t,\cdot)>\e)\}}(x,t)\,dxdt \\
{}+{}\int\int \left(\int v(f_h-\Pi_{f_h})\,dv\right)w(x,t)\id_{\{\int \Pi_{f_h}\,dv\geq h^{1/4}\}\cap\{D(f_h(x,t,\cdot)\geq\e)\}}(x,t)\,dxdt.
\end{multline*}
The first integral is of the order $h^{1/4},$ while the second is of order $h^{1/2}$ because of \eqref{est:2.2} and 
\[
\int v\Pi_{f_h}\,dv{}={}\left(\int \Pi_{f_h}\,dv\right)^2\geq h^{1/2}.
\]

\end{section}

\begin{section}{An example}
As an example of the approximate solutions, consider one-dimensional Burger's equation
\[
\partial_t u{}+{}\partial_x(u^2/2){}={}0,
\]
with initial data
\[
u_0{}={}
\left\{
\begin{array}{rl}
1, & x<0,\\
-1, & x>0.
\end{array}
\right.
\]
The admissible weak solutions to this problem is a stationary shock wave: 
\begin{equation}
\label{eq:shock}
\forall t>0,\, u_s(x,t){}={}
\left\{
\begin{array}{rl}
1, & x<0,\\
-1, & x>0.
\end{array}
\right.
\end{equation}

We're going to compute the time-discrete solution by the given by  \eqref{step:transport}, \eqref{step:collapse}. Let $\e>0.$ For the collapse step, we will need only to compute the equilibrium kinetic density for a density in one of  the forms
\[
f(v){}={}-\id_{[-1,0]}(v){}+{}\id_{[a,1]}(v),\,
f(v){}={}\id_{[0,1]}(v){}-{}\id_{[-a,0]}(v),\,a>0.
\]
The corresponding equilibrium kinetic functions are
\[
f_0(v){}={}-\id_{[-1+a,0]}(v),\, \id_{[1-a,0]}(v).
\]
From the condition \eqref{step:small} we obtain that when $\e$ is small, there is $\delta{}={}O(\e)\in(0,1),$ such that $f(v)$ is replaced at the collapse step  with $f_0(v)$ if and only if $a<1-\delta.$ 

To make computations a bit easier we take the initial kinetic density to be
\[
f_h(x,0,v){}={}
\left\{
\begin{array}{rl}
1, & v\in[0,1],\,x\leq \max\{0, (v-(1-\delta))h\},\\
0, & v>1,\,x\in\mathbb{R},\\
0, & v<-1,\,x\in\mathbb{R},\\
-1, & v\in[-1,0],\, x\geq\min\{0, (v+(1-\delta))h\}.
\end{array}
\right.
\]
The support of $f_h(x,0,v)$ is sketched on figure 3.1

Function $u_{0,h}(x)=\int f_h(x,0,v)\,dv$ equals $u_0(x)$ outside of small interval $(-h,h),$ and converges to it  pointwisely and in any $L^p_{loc}(\mathbb{R})$ norm. The solutions of Burger's equation with initial data $u_{0,h}$ converge to the shock \eqref{eq:shock} as $h\to 0.$

Denote
\[
A_0^+{}={}\{(x,v)\,:\, f(x,0,v){}={}1\},\, A_0^-{}={}\{(x,v)\,:\, f(x,0,v){}={}-1\}.
\]
Let $T^h\;:\;\mathbb{R}^{2}_{x,v}\to \mathbb{R}^{2}_{x,v}$
be the translation in $x$ with velocity $v:$
\[
T^t(x,v){}={}(x+tv,v).
\]
For time $t\in(0,h),$ 
\[
f_h(x,t,v){}={}\id_{T^t[A_0^+]}(x,v){}-{}\id_{T^t[A_0^-]}(x,v).
\]
After the collapse step at $t=h,$ density $f_h$ equals $1$ on the set
\[
A_1^+{}={}\{ (x,v)\;:\; v\in[0,1],\,x\leq\max\{ -h(1-\delta),-hv\}\} \cup \Delta^+,
\]
\[
\Delta^+{}={}\{ (x,v)\;:\; v\in[1-\delta,1],\, (1-\delta)h\leq x\leq 2hv\},
\]
and equals $-1$ on the set
\[
A_1^-{}={}\{ (x,v)\;:\; v\in[-1,0],\,x\geq\min\{ h(1-\delta),-hv\}\} \cup \Delta^-,
\]
\[
\Delta^-{}={}\{ (x,v)\;:\; v\in[-1,-(1-\delta)],\, -2hv\leq x\leq -(1-\delta)h\}.
\]
Thus, we can write
\[
f_h(x,h,v){}={}\id_{A_1^+}(x,v)-\id_{A_1^-}(x,v)+\id_{T^0[\Delta^+]}(x,v){}-{}\id_{T^0[\Delta^-]}(x,v).
\]
The support of $f_h$ is sketched on figures 3.2.

Then, for $t\in[h,2h):$
\[
f_h(x,t,v){}={}\id_{T^{t-h}[A_1^+]}(x,v)-\id_{T^{t-h}[A_1^-]}(x,v)+\id_{T^{t-h}[\Delta^+]}(x,v){}-{}\id_{T^{t-h}[\Delta^-]}(x,v).
\]
At, $t=2h,$ after the interaction step,
\[
f_h(x,2h,v){}={}\id_{A_0^+}(x,v)-\id_{A_0^-}(x,v)+\id_{T^h[\Delta^+]}(x,v){}-{}\id_{T^h[\Delta^-]}(x,v),
\]
so that for $t\in[2h,3h),$
\[
f_h(x,t,v){}={}\id_{T^{t-2h}[A_0^+]}(x,v)-\id_{T^{t-2h}[A_0^-]}(x,v)+\id_{T^{t-h}[\Delta^+]}(x,v){}-{}\id_{T^{t-h}[\Delta^-]}(x,v),
\]
and at $t=3h,$
\begin{multline*}
f_h(x,3h,v){}={}\id_{A_1^+}(x,v){}-{}\id_{A_1^-}(x,v){}+{}\id_{T^{2h}[\Delta^+]}(x,v){}-{}\id_{T^{2h}[\Delta^-]}(x,v)\\
+\id_{T^0[\Delta^+]}(x,v){}-{}\id_{T^0[\Delta^-]}(x,v).
\end{multline*}
The support of $f_h$ at time $t=2h$ is sketched on figure 3.3.

If follows that for $t\in[(2k-1)h,2kh),$
\begin{multline*}
f_h(x,t,v){}={}\id_{T^{t-(2k-1)h}[A_1^+]}(x,v){}-{}
\id_{T^{t-(2k-1)h}[A_1^-]}(x,v)\\
+\sum_{i=1}^k \id_{T^{t-(2i-1)h}[\Delta^+]}(x,v)
-\id_{T^{t-(2i-1)h}[\Delta^-]}(x,v),
\end{multline*}
and for $t\in[2kh,(2k+1)h),$
\begin{equation*}
f_h(x,t,v){}={}\id_{T^{t-2kh}[A_0^+]}(x,v){}-{}
\id_{T^{t-2kh}[A_0^-]}(x,v)
+\sum_{i=1}^k \id_{T^{t-(2i-1)h}[\Delta^+]}(x,v)
-\id_{T^{t-(2i-1)h}[\Delta^-]}(x,v).
\end{equation*}
Form these formulas we obtain the following information. When $v\in[1-\delta,1]$ and $x<0,$ then $f_h{}={}1;$ when
$v\in[1-\delta,1]$ and $x>vt,$  $f_h{}={}0;$ for $v\in[1-\delta,1]$ and $x\in(h,vt),$ $f_h$ (as a function of $x$) has a form of a square wave: it equals $1$ on intervals of length $2h(v-(1-\delta))$ separated by intervals of legnth $2hv,$ where function equals $0.$ 

When $v\in[-1,-(1-\delta)]$ and $x>0,$ then $f_h{}={}-1;$ when
$v\in[-1,-(1-\delta)]$ and $x<vt,$  $f_h{}={}0;$ for $v\in[-1,-(1-\delta)]$ and $x\in(vt,-h),$ $f_h$ (as a function of $x$) has a form of a square wave: it equals $-1$ on intervals of length $2h(v-(1-\delta))$ separated by intervals of legnth $2hv,$ where function equals $0.$ 

For $v\in(-(1-\delta),1-\delta),$ $f_h$ equals either
\[
\id_{T^{t-(2k-1)h}[A_1^+]}(x,v){}-{}\id_{T^{t-(2k-1)h}[A_1^-]}(x,v),
\]
or
\[
\id_{T^{t-2kh}[A_0^+]}(x,v){}-{}\id_{T^{t-2kh}[A_0^-]}(x,v).
\]

The following types of convergence take place when $h\to0:$
\begin{enumerate}
\item $\forall v\in(-(1-\delta),1-\delta),\,t>0:$
\[
f_h(x,t,v){}\to
\left\{\begin{array}{rl}
1, & x<0,\,v\in[0,1-\delta),\\
0, & x>0,\,v\in[0,1-\delta),\\
0, & x<0,\, v\in(-(1-\delta),0],\\
-1, & x>0,\,v\in(-(1-\delta),0];
\end{array}
\right.
\]
\item $\forall v\in[(1-\delta),1],\,t>0:$
 \[
f_h(\cdot,t,v){}\to
\left\{\begin{array}{ll}
1, & x<0,\\
\dfrac{v-(1-\delta)}{2v-(1-\delta)}, & 0<x<vt,\\
0, & x<vt,
\end{array}
\right.
\]
in *--weak tolplogy of $L^\infty_{loc}(\mathbb{R}_x);$
\item $\forall v\in[-1,-(1-\delta)],\,t>0:$
 \[
f_h(\cdot,t,v){}\to
\left\{\begin{array}{ll}
0, & x<vt,\\
\dfrac{v+(1-\delta)}{2v+(1-\delta)}, & vt<x<0,\\
-1, & x>0,
\end{array}
\right.
\]
in *--weak topology of $L^\infty_{loc}(\mathbb{R}_x).$
\end{enumerate}
The limiting function $f$ has self-similar structure, depending on $x/t.$ For $x>0$ it is given by the formulas
\[
f(x,t,v){}={}\left\{
\begin{array}{ll}
-\id_{[-1,0]}(v), & x/t\geq 1,\\
{}-{}\id_{[-1,0]}(v)+\dfrac{v-(1-\delta)}{2v-(1-\delta)}\id_{[x/t,1]}(v), & 1-\delta<x/t<1,\\
{}-{}\id_{[-1,0]}(v)+\dfrac{v-(1-\delta)}{2v-(1-\delta)}\id_{[1-\delta,1]}(v), & 0<x/t<1-\delta,
\end{array}
\right.
\]
while for $x<0$ by 
\[
f(x,t,v){}={}\left\{
\begin{array}{ll}
\id_{[0,1]}(v), & x/t\leq -1,\\
\id_{[0,1]}(v)+\dfrac{v+(1-\delta)}{2v+(1-\delta)}\id_{[-1,x/t]}(v), & -1<x/t<-(1-\delta),\\
\id_{[0,1]}(v)+\dfrac{v+(1-\delta)}{2v+(1-\delta)}\id_{[-1,-(1-\delta)]}(v), & -(1-\delta)<x/t<0.
\end{array}
\right.
\]

This should be compared to the kinetic density of the shock \eqref{eq:shock}:
\[
f_s(x,t,v){}={}\left\{
\begin{array}{cl}
-\id_{[-1,0]}(v), & x/t>0,\\
\id_{[0,1]}(v), & x/t<0.
\end{array}
\right.
\]
\end{section}
Notice, that there is positive mass leaking through the shock to the right of $x=0$ at speeds close to $v=1,$ and,  negative mass leaking to the left of  $x=0.$

\begin{figure}
\begin{center}
\includegraphics[trim=0 450 0 0, scale=.75]{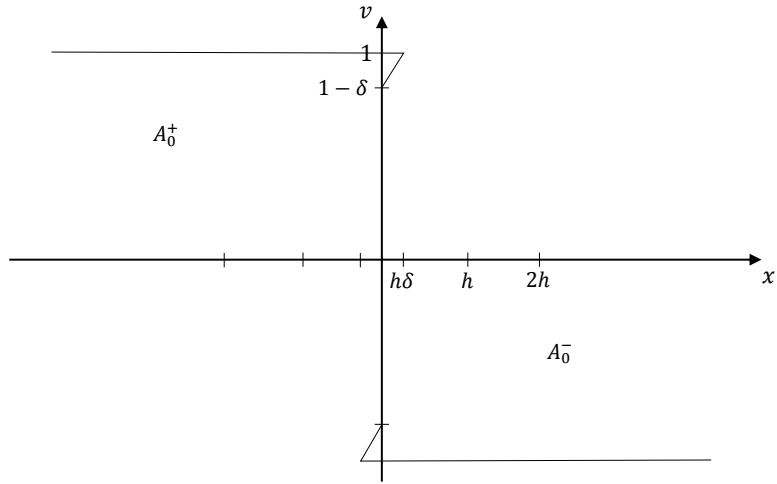}        
\caption{The support of the kinetic density at time $t=0.$}
\end{center}
\label{fig1}
\end{figure}

\begin{figure}
\label{fig2}
\begin{center}
\includegraphics[trim=0 450 0 0, scale=.75]{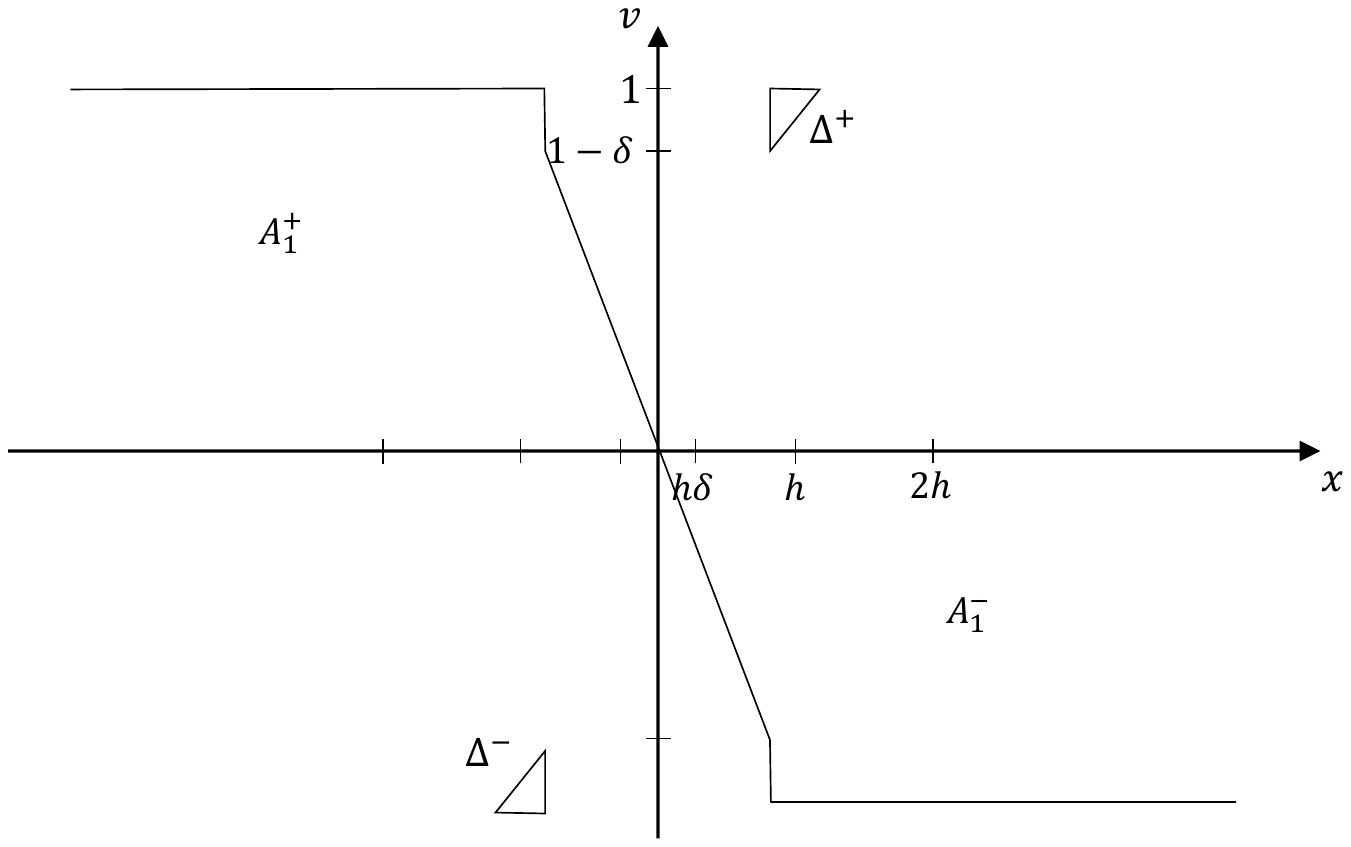}        
\caption{The support of the kinetic density at time $t=h.$}
\end{center}
\end{figure}

\begin{figure}
\label{fig3}
\begin{center}
\includegraphics[trim=0 450 0 0, scale=.75]{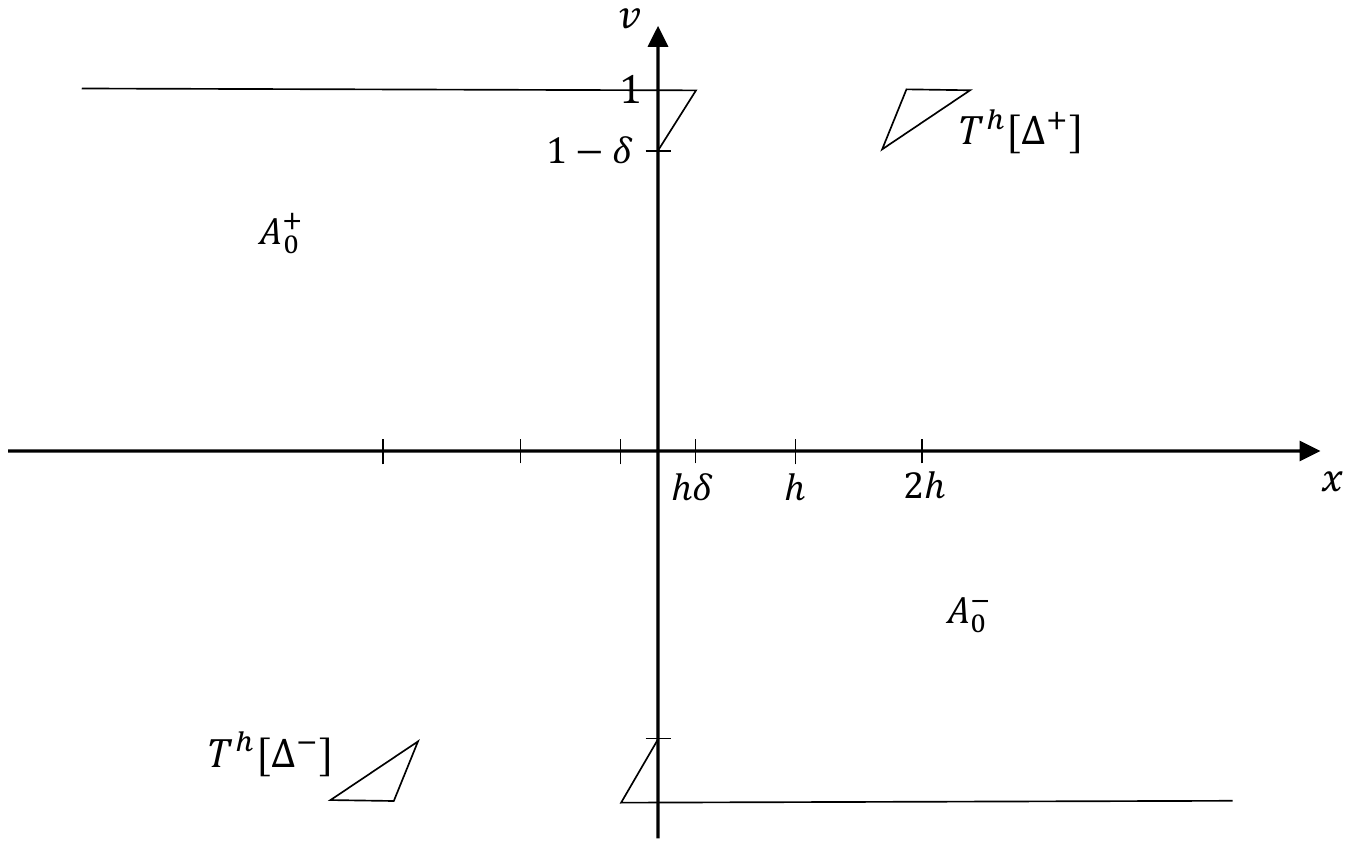}        
\caption{The support of the kinetic density at time $t=2h.$}
\end{center}
\end{figure}

\begin{section}{Appendix}
\begin{lem}
\label{lemma:ap1}
Let $f\in L^1((0,+\infty)),$ with bounded support and such that $f(v)\in[0,1],$ for a.e. $v.$ Let $f_0(v){}={}\Pi_{f(\cdot)}(v).$
There is a non-decreasing map
\[
\tau_{f_0}^f\;:\;[0,+\infty)\to[0,+\infty),
\]
such that
\[
\tau_{f_0}^f(v)\geq v,\quad \forall v\geq0,
\]
and for all test function $\phi\in L^\infty(\mathbb{R}^+_v),$
\begin{equation}
\label{lemma_ap:1}
\int \phi(\tau_{f_0}^f(v))f_0(v)\,dv{}={}\int \phi(v) f(v)\,dv.
\end{equation}
In particular, if $\eta$ is a convex function, then
\begin{equation}
\label{lemma_ap:2}
\int \eta'(v) f(v)\,dv\geq \int \eta'(v) f_0(v)\,dv.
\end{equation}
\end{lem}
\proof
The map $\tau_{f_0}^f$ in the statement of the lemma can be taken to be the optimal transport map from the measure $f_0dv$ to $fdv.$ It is constructed in the following way. Set
\[
F(v){}={}\int_0^v f(v)\,dv,\quad F_0(v){}={}\int_0^v f_0(v)\,dv.
\]
The map $\tau_f^{f_0}(v){}={}F(v)$ is a mass transport map from $f\,dv$ to $f_0\,dv:$ for any $\phi\in L^\infty(\mathbb{R}^+_v),$
\[
\int \phi(\tau_f^{f_0}(v))f(v)\,dv{}={}\int \phi(v) f_0(v)\,dv.
\]
Moreover, $\forall v\geq0,$ $\tau_f^{f_0}(v)\leq v.$
The generalized inverse of $\tau_f^{f_0},$
\[
\tau_{f_0}^f(v){}={}\inf\{ u\,:\, F(u)>v\}
\]
is a  mass transport map from $f_0\,dv$ to $f\,dv$ and
\[
\tau_{f_0}^f(v)\geq v,\quad \forall v\geq 0.
\]
Inequality \eqref{lemma_ap:2} follows from \eqref{lemma_ap:1}

A variant of the following lemma was used in \cite{Brenier}.
\begin{lem}
\label{lemma:Brenier}
Given $\eta(v)$ -- a convex function on $\mathbb{R}_v$ and two bounded measurable sets
\[
A_1\subset [0,+\infty),\quad A_2\subset(-\infty,0],
\]
it holds
\[
 \int \eta'(v)\id_{{\rm Int}(0,|A_1|-|A_2|)}(v)\,dv{}\leq{}\int \eta'(v)(\id_{A_1}(v)-\id_{A_2}(v))\,dv.
\]
\end{lem}

\proof
Let $\eta$ be a convex function on $\mathbb{R}_v.$ By subtracting from it a supporting line to the epigraph of $\eta$ at point $v=0,$ we may assume that $\eta'(v)\geq0$ for $v\geq 0$ and $\eta'(v)\leq 0,$ for $v\leq 0$ and $\eta(0){}={}0.$

Suppose first, that $A_2$ is empty. The statement of the lemma follows then from \eqref{lemma_ap:2}. When $A_2$ is not empty, by symmetry,
\[
-\int \eta'(v)\id_{A_2}(v)\,dv{}\geq{}-\int \eta'(v)\id_{[-|A_2|,0]}\,dv.
\]
It suffices to consider the case $|A_1|>|A_2|.$ We estimate:
\begin{multline*}
\int \eta'(v)(\id_{A_1}(v)-\id_{A_2}(v))\,dv{}\geq{}
\int \eta'(v)\id_{[0,|A_1|]}\,dv{}-{}\int \eta'(v)\id_{[-|A_2|,0]}\,dv\\
{}={}\eta(|A_1|){}+{}\eta(-|A_2|)\geq \eta(|A_1|-|A_2|){}={}\int \eta'(v)\id_{[0,|A_1|-|A_2|]}(v)\,dv.
\end{multline*}

The following theorem was proved in \cite{BGP}, following \cite{Gerard}.
\begin{thm}
\label{Gerard}
Let $\Omega$ be an open set of $\mathbb{R}^+_{x,t},$ $a\in L^\infty_{loc}(\mathbb{R}_v)^d$ such that 
\begin{equation}
\label{nonlinear}
\forall\, \sigma \in \mathbb{S}^{d-1}\quad \forall \xi\in\mathbb{R}\quad |\{ v\in\mathbb{R}\;:\; a(v)\cdot\sigma{}={}\xi\}|{}={}0,
\end{equation}
and $f_n$ a bounded sequence in $L^2_{loc}(\Omega\times\mathbb{R}_v)$ such that
\begin{equation}
\label{Gerard:1}
\forall \psi\in C^\infty_c(\mathbb{R}_v)\quad \int (\partial_t f_n{}+{}a(v)\cdot\grad_x f_n)\psi\,dv\quad \mbox{is pre-compact in } H^{-1}_{loc}(\Omega).
\end{equation}
Then,
\[
\forall \psi\in C^\infty_c(\mathbb{R}_v)\quad \int f_n\psi\,dv\quad\mbox{is pre-compact in } L^2_{loc}(\Omega).
\]
\end{thm}

Weak solutions of the problem
\begin{equation}
\label{eq:CBGK}
\partial_t f{}+{}A'(v)\cdot\grad_x f{}={}(\Pi_{f(x,t,\cdot)}(v)-f(x,t,v))\id_{\{D(f(x,t,\cdot)>\e\}}(x,t),\quad \mathcal{D}'(\mathbb{R}^+_{x,t,v}),
\end{equation}
\begin{equation}
\label{eq:CBGK_initial}
f(x,0,v){}={}f_0(x,v).
\end{equation}

Solution of this problem is obtained by approximating the right-hand side of \eqref{eq:CBGK}
by a smoother operator and using the theory of $C^0$-- semigroups. 

Let $r>0$ be an approximating parameter. Let $\id^r(s)$ be a smooth approximation of the Heavyside function: $\id^r(s){}={}0, $ for $s\leq0;$ $\id^r(s){}={}1$, for $s\geq 2r,$ and $\id^r(s)$ is continuously differentiable and monotone. 

Let $w$ be a cut-off function
\[
w(s){}={}s\id_{[0,1]}(s){}+{}\id_{\{s>1\}}(s).
\]

Let $X=L^1(\mathbb{R}^d\times[0,M]),$ for some $M>0.$ 
Define an operator $L^r\;:\;X\to X,$ by
\[
L^r(f){}={}(\Pi_{w\circ f(x,t,\cdot)}-w\circ f)\id^r\circ\left(
\int v( w\circ f(x,t,v)-\Pi_{w\circ f(x,t,\cdot)}(v))\,dv-\e\int v\Pi_{w\circ f(x,t,\cdot)}(v)\,dv\right).
\]

One can easily verify the next lemma.
\begin{lem} Operator $L^r$ is uniformly Lipschitz on $X:$ there is $C>0,$ such that for any $f,g\in X,$ \[
\|L^r(f)-L^r(g)\|_X{}\leq{} C\|f-g\|_X.
\]
\end{lem}

Consider a semigroup $S_t$ on $X$ is defined by
\begin{equation}
\label{def:St}
S_t(f)(x,v){}={}f(x-A'(v) v,v).
\end{equation}

\begin{lem} Let $A$ be Lipschitz continuous on $[0,M].$ $S_t$ is a $C^0$--semigroup with the generator  $A'(v)\cdot\grad_x$ with the domain 
\[
D{}={}\{f\in X\;:\; A'(v)\cdot\grad_x f\in X\}.
\]
\end{lem}
For the proof of this lemma, see for example proposition 1, ch. XXI, of \cite{LD}.
Consider the problem
\begin{equation}
\label{eq:SG}
\partial_t f{}+{}A'(v)\cdot\grad_x f{}={}L^r(f),\quad \mbox{in }\mathcal{D}'(\mathbb{R}^+_{x,t}),
\end{equation}

\begin{equation}
\label{eq:SG_initial}
f(x,0,v){}={}f_0(x,v).
\end{equation}

\begin{thm}
\label{th:}
Let $f_0\in X,$ such that a.e. $(x,v),$
\[
f_0(x,v)\in[0,1].
\]
There is exist a unique weak solution $f\in C([0,+\infty); X)$ of problem \eqref{eq:SG}, \eqref{eq:SG_initial}.
For all $t>0$ a.e. $(x,v),$
\[
f(x,t,v)\in[0,1],
\]
and 
\[
\sup_{t>0}\|f(\cdot,t,\cdot)\|_{X}\leq \|f_0(\cdot)\|_X.
\]
If for some $M_0\in(0,M),$ and a.e. $x,$
\[
\supp f_0(x,\cdot)\subset [0,M_0],
\]
then for all $t>0$ and a.e. $(x,v),$
\[
\supp f(x,t,v)\subset [0,M_0].
\]
\end{thm}

\proof
Using theorem 1.2, ch. 6 of \cite{Pazy} we obtain a unique mild solution $ f\in C([0,+\infty); X):$
for any $t>0$
\[
f(t,\cdot){}={}S_t(f_0(\cdot)){}+{}\int_0^t S_{t-s}\circ L^r(f(\cdot,s,\cdot))\,ds.
\]
Using definition \eqref{def:St} of $S_t$ we conclude that $f$ is a weak solution (solution in a sense of distributions) of \eqref{eq:SG}. Since the transport velocity is independent of $x$ and $f$ and $L^r(f)$ are uniformly bounded, $f$ is a renormalized weak solution of \eqref{eq:SG}, see \cite{DiPerna-Lions}: if $w\in C^1(\mathbb{R}),$ then $w\circ f$ is a weak solution of 
\[
\partial_t(w\circ f){}+{}A'(v)\cdot\grad_x (w\circ f){}={}L^r(f) (w'\circ f).
\]
Selecting $w(s){}={}(s-1)^2\id_{\{s>1\}}(s),$  and using the fact that if $f(x,t,v)>1, $ then $L^r(f)(x,t,v)\leq0,$ we find after integration that for all $t>0,$
\[
\int\int_{0}^M (f-1)^2\id_{\{f>1\}}\,dxdv\leq 0,
\]
so that $f\leq 1$ a.e. $(x,v).$ A similar argument shows that $f\geq0$ a.e. $(x,v).$

The statement about the support of $f(x,t,\cdot),$ follows from the fact that if $\supp f(x,\cdot)\subset[0,M_0]$ for a. e. $x,$ then $\supp L^r(f)(x,\cdot)\subset [0,M_0],$ a.e. $x.$

The estimate on $\sup_t\|f(x,t,v)\|_X$ follows since $\int L^r(f(x,t,v))\,dv{}={}0,$ a.e. $x,$ and $f\geq0.$

When the range of the solution from the last theorem $f\in[0,1],$ $f$ is a renormalized weak solution of equation without the cut-off function $w:$
\begin{equation}
\partial_t f{}+{}A'(v)\cdot\grad_x f{}={}(\Pi_f-f)\id^r\circ\left( \int v(f-\Pi_f)\,dv -\e\int v\Pi_{w\circ f(x,t,\cdot)}(v)\,dv\right).
\end{equation}
The solution $f$ depends on two parameters $r$ and $M.$ We set $M=1/r$ and take limit
for the sequence $\{f^r\}$ as $r\to0,$ to obtain a weak solution of \eqref{eq:CBGK}, \eqref{eq:CBGK_initial}.
\begin{thm}
\label{th:last}
Let $A$ verify non-degeneracy condition \eqref{main:nondegeneracy}. Let $f_0\in L^1(\mathbb{R}^{d+1}_{x,v}),$
such that $f_0\in[0,1],$ a. e. $(x,v),$ and for a.e. $x,$
\[
\supp f_0(x,\cdot)\subset [0,M_0]
\] 
for some $M_0$ independent of $x.$
There exists a weak solution 
\[ f\in L^\infty(\mathbb{R}^+_{x,t,v})\cap L^\infty(0,+\infty;L^1(\mathbb{R}^{d+1}_{x,v}))\cap C([0,+\infty); L^p_{loc, weak}(\mathbb{R}^{d+1}_{x,v})),\]
 $1\leq p<+\infty,$  of \eqref{eq:CBGK}, \eqref{eq:CBGK_initial}.
For all $t>0,$ $f(x,t,v)\in[0,1],$ a.e. $(x,v)$ and for a.e. $x,$
\[
\supp f(x,t,\cdot)\subset [0,M_0].
\]
\end{thm}
\proof Consider a family of solutions $\{f^r\}$ of problem \eqref{eq:SG}, \eqref{eq:SG_initial} with initial data $f_0.$ $f^r$ is bounded in $L^\infty(\mathbb{R}^+_{x,t,v})\cap L^\infty(0,+\infty;L^1(\mathbb{R}^{d+1}_{x,v})),$ uniformly in $r.$ It follows that there is 
\[
f\in L^\infty(\mathbb{R}^+_{x,t,v})\cap L^\infty(0,+\infty;L^1(\mathbb{R}^{d+1}_{x,v}))
\]
 and 
a sequence (still labeled by $r$) such that $f^r\to f$ in *--weak topologies of $L^\infty(\mathbb{R}^+_{x,t,v})$ and $L^\infty(0,+\infty;L^1(\mathbb{R}^{d+1}_{x,v})).$
The right-hand sides $L^r$ of equation \eqref{eq:SG} are bounded in $L^\infty(\mathbb{R}^+_{x,t}).$
Using theorem \ref{Gerard}, we conclude that for any $\psi\in C(\mathbb{R}_v),$
\[
\left\{ \int \psi f^r\,dv\right\} \quad \mbox{pre-compact in } L^p_{loc}(\mathbb{R}^+_{x,t}),
\]
for any $1\leq p<+\infty.$ In particular, there is a further subsequence, such that
\[
\int v f^r\,dv\to \int v f\,dv\quad \mbox{a.e. $(x,t)$ and in  $L^p_{loc}(\mathbb{R}^+_{x,t})$}
\]
and $\Pi_{f^r}\to \Pi_f$ a.e. $(x,t,v)$ and  any $L^p_{loc}(\mathbb{R}^+_{x,t,v}),$ $1\leq p<+\infty.$
It follows that
\[
(\Pi_{f^r}-f^r)\id^r\circ\left( \int v(f^r-\Pi_{f^r})dv-\e\int v\Pi_{f^r}\,dv\right)
 \to (\Pi_f -f)\id^0\left( \int v(f-\Pi_{f})dv-\e\int v\Pi_f\,dv\right)
\]
in *--weak $L^\infty_{loc}(\mathbb{R}^+_{x,t,v}).$ Thus, $f$ is a weak solution of \eqref{eq:CBGK}.

\end{section}


\begin{thebibliography}{XXX}


\bibitem{BB} F.~Berthlin and F.~Bouchut, {\sl Relaxation to isentropic gas dynamics for a BGK system with
single kinetic entropy,\/}
 Methods Appl. Anal., 9 (2002), p. 313--327.


\bibitem{BGP} F.~Bouchut, F.~Golse and M.~Pulvirenti, {\sl Kinetic equations and asymptotic theory,\/} Series in Applied Math. Series Editors  P.G.~Ciarlet and P.-L.~Lions, Gauthier-Villars (1998).

\bibitem{Brenier} Y.~Brenier, {\sl Averaged multivalued solutions for scalar conservation laws,\/} SIAM J. Numer. Anal. 27 (1984), n. 6, p. 1013--1037.



\bibitem{LD} R.~Dautray and J.-L.~Lions, {\sl Mathematical Analysis and Numerical Methods for Science and Technology,\/} vol. 6, Springer-Verlag, (1993).

\bibitem{Constantine} C.~M.~Dafermos, {\sl Hyperbolic Conservation Laws in Continuum Physics,\/} Springer, (2010).


\bibitem{DiPerna} R.~DiPerna, {\sl Convergence of the viscosity method for isentropic gas dynamics,\/} Comm. Math. Phys. 91 (1983), p. 1--30.

\bibitem{DiPerna-Lions} R.~DiPerna, P.L.~Lions, {\sl Ordinary differential equations, transport theory and Sobolev spaces,\/} Invent. Math.  98, (1989), 3, p. 511--547.


\bibitem{Gerard} P.~G\'{e}rard, {\sl Microlocal defect measures,\/} Comm. Partial Diff. Eq.  16 (1991), p. 1761--1794.

\bibitem{GM} Y.~Giga, R.~Miyakawa, {\sl A kinetic construction of global solutions of first order quasilinear equations,\/} Duke Math. J., 50 (1983), p. 505--515.


\bibitem{GLPS} F.~Golse, P.-L.~Lions, B.~Perthame and R.~Sentis, {\sl Regularity of the moments of the solution of a transport equation,\/} J. Funct. Anal. 76 (1988), p.110--125.


\bibitem{LPT} P.-L.~Lions, B.~Perthame, E.~Tadmor, {\sl Kinetic formulation of the multi-dimensional scalar conservation laws,\/} J. of AMS 7 (1994), no. 1, p.196--191.

\bibitem{LPT2} P.-L.~Lions, B.~Perthame, E.~Tadmor, {\sl Kinetic formulation of the isentropic gas dynamics and p
p
-systems,\/} Comm. Math. Phys. 163 (1994), no. 2, p. 415--431.

\bibitem{Pazy} A.~Pazy, {\sl Semigroups of Linear Operators and Applications to Partial Differential Equations,\/}
 Applied Math. Sciences vol. 44, Springer (1983).

\bibitem{Perepelitsa_isentropic}M.~Perepelitsa, {\sl A kinetic model for the approximately isentropic solutions of the Euler equations,\/} to appear in JDE,  preprint available at www.math.uh.edu/~misha, (2015).

\bibitem{Perthame} B.~Perthame, {\sl Global Existence to the BGK Model of Boltzmann Equation,\/} J. Diff. Eqations 82 (1989), p. 191--205.



\bibitem{Perthame_m} B.~Perthame, {\sl Kinetic Formulation of Conservation Laws,\/} Oxford Lecture Ser. Math. Appl.
21, Oxford University Press, New York, (2002).


\end{thebibliography}
\end{document}